\documentclass[oneside,12pt]{amsart}
\usepackage{amssymb, amscd}
\usepackage[all]{xy}

\newtheorem{thm}{Theorem}[section]
\newtheorem{cor}[thm]{Corollary}
\newtheorem{lem}[thm]{Lemma}
\newtheorem{defn}[thm]{Definition}
\newtheorem{prop}[thm]{Proposition}
\newtheorem{remk} [thm]{Remark}
\newtheorem{examp}[thm]{Example}

\newtheorem{conj}[thm]{Conjecture}

\newcommand{\del}[2]{{}}

\title{Jumps in Cohomology and Free Group Actions}

\author{Nansen Petrosyan}
\address{Department of Mathematics, Indiana University, Bloomington IN 47405, USA}%
\email{nanpetro@indiana.edu}%

\thanks{}%

\date{\today}

\begin{document}
\begin{abstract}
{\footnotesize A discrete group $G$ has periodic cohomology over $R$
if there is an element in a cohomology group, cup product with which
induces an isomorphism in cohomology after a certain dimension. Adem
and Smith showed if $R={\mathbb Z}$, then this condition is
equivalent to the existence of a finite dimensional
free-$G$-CW-complex homotopy equivalent to a sphere. It has been
conjectured by Olympia Talelli, that if $G$ is also torsion-free
then it must have finite cohomological dimension. In this paper we
use the implied condition of jump cohomology over $R$ to prove the
conjecture for ${\sl H}{\mathcal F}$-groups and solvable groups. We
also find necessary conditions for free and proper group actions on
finite dimensional complexes homotopy equivalent to closed,
orientable manifolds.}
\end{abstract}

\maketitle
\section{Introduction}

\vspace{.2cm}

\begin{defn}\label{LHS} A discrete group $G$ has {\it periodic cohomology over} $R$ if
there is a cohomology class $\alpha\in$ Ext$_{RG}^{\ast}(R,R)$ with
$|\alpha|>0$ and an integer $n\geq 0$ such that the cup product map
\begin{center}
 $\alpha \cup -:$ Ext$_{RG}^{i}(R, M)\rightarrow$ Ext$_{RG}^{i+|\alpha|}(R, M)$
\end{center}
\noindent is an isomorphism for every $RG$-module $M$ and every
integer $i\geq n$.
\end{defn}

If $R={\mathbb Z}$, then we simply say that $G$ has {\it periodic
cohomology}. In this case Adem and Smith \cite{AS} proved that this
condition is equivalent to the existence of a finite dimensional
free-$G$-CW-complex homotopy equivalent to a sphere. If a group $G$
has periodic cohomology then it has periodic cohomology over any
commutative ring $R$ with a unit. This definition of periodicity
using the cup product differs from the general definition of periodic
cohomology of a group. Namely, a group $G$ has {\it periodic
cohomology after $k$-steps} if there exists an integer $q\geq 0$
such that the functors  $H^{i}(G, -)$ and $H^{i+q}(G, -)$ are
naturally equivalent for all $i>k$.  We will call this classical
periodic cohomology. It is immediate that the cup product notion of
periodic cohomology implies the classical periodic cohomology. It is
a conjecture by O. Talelli that the two notions are the same. She
has proved that the conjecture holds in the case of the ${\sl
H}{\mathcal F}$-groups introduced by P. Kropholler in \cite{krop}.
This is the smallest class of groups which contains all finite
groups and which contains all groups $G$, whenever $G$ acts
cellularly on a finite dimensional contractible CW-complex with all
isotropy subgroups already in ${\sl H}{\mathcal F}$. This is a large
class of groups. Among others, it contains all countable linear
groups, all countable solvable groups and all groups with finite
virtual cohomological dimension.

The following was conjectured by Talelli in \cite{tal2}, Conjecture
III:

\vspace{.2cm}

\noindent {\bf Conjecture} (O. Talelli). Every torsion-free discrete
group $G$ with classical periodic cohomology after some steps has
finite cohomological dimension. \vspace{.2cm}

The Conjecture was shown to be true for the class of ${\sl
H}{\mathcal F}$-groups by Mislin and Talelli in \cite{MT}. In our
paper we investigate a weaker condition, implied by cohomological
periodicity, for a discrete group $G$.

\vspace{.2cm}
\begin{defn}\label{LHS} A discrete group $G$ has {\it jump cohomology over} $R$ if there
exists an integer $k\geq 0$, such that if $H$ is any subgroup of $G$
with $cd_R(H)<\infty$, then $cd_R(H)\leq k$. The bound $k$ will be
called the {\it jump height over R}. If $R={\mathbb Z}$ then we say
$G$ has {\it jump cohomology} and a {\it jump height k}.
\end{defn}

\vspace{.2cm}

By naturality of the cup product it follows that if a group $G$ has
periodic cohomology over $R$, then $G$ has jump cohomology over $R$.
In section 2 we prove the following proposition which relates the
geometric aspects of the notion to the jump cohomology.

\vspace{.2cm}

\begin{prop}\label{LHS} Let $G$ be a discrete group that acts freely and
properly on an
$e$-dimensional CW-complex $E$. Suppose there is an integer $n\leq e$
such that $H_i(E,{\mathbb Z})=0$ for all $i>n$, and
$H_n(E,{\mathbb Z})\cong {\mathbb Z}$. Set $k=e-n$. If $H$ is any
subgroup of $G$ with $cd(H)< \infty$ ($hd(H)< \infty$),
then $cd(H)\leq k$ (respectively $hd(H)\leq k$).
\end{prop}

\vspace{.2cm}

In section 3 we investigate the classes of ${\sl
H}\mathcal{F}$-groups  and solvable groups without $R$-torsion.  We
prove the following two theorems. \vspace{.2cm}
\begin{thm}\label{LHS} Let $G$ be a discrete group without $R$-torsion
with jump cohomology of height
$k$ over $R$. If $G$ is an ${\sl H}\mathcal{F}$-group, then
$cd_R(G)\leq k$. In particular, any ${\sl H}\mathcal{F}$-group $G$
has jump cohomology of height $k$ over $\mathbb Q$ if and only if
$cd_{\mathbb Q}(G)\leq k$.
\end{thm}

\vspace{.2cm}

\begin{thm}\label{LHS} Let $G$ be a solvable group that acts freely
and properly
on an $e$-dimensional CW-complex $E$. Suppose there is an integer $n\leq e$
such that $H_i(E,{\mathbb Z})=0$ for all $i>n$, and $H_n(E,{\mathbb Z})\cong
{\mathbb Z}$, then $h(G)\leq e - n$. If $G$ is also $R$-torsion-free,
$hd_R(G)\leq e - n$.
\end{thm}

\vspace{.2cm}

Our results give affirmative answers for the case of ${\sl
H}\mathcal{F}$-groups, and when $R$ is a domain of characteristic
zero for the case of solvable groups to the conjecture we believe is
a natural generalization of the conjecture by Talelli.

\vspace{.2cm}

\begin{conj}\label{LHS} For every discrete group $G$ without $R$-torsion the following are
equivalent.

\begin{enumerate}

\item $G$ has jump cohomology of height $k$ over $R$.

\vspace{1mm}
\item  $G$ has periodic cohomology over $R$ starting in dimension $k+1$.

\vspace{1mm}
\item $cd_R(G)\leq k$
\end{enumerate}
\end{conj}

This conjecture in particular would imply Talelli's Conjecture. It
also says that every group with periodic cohomology over ${\mathbb
Q}$ must have $cd_{\mathbb Q}(G)<\infty$. Periodic or jump
cohomology is a rather strong condition for the cohomology of $G$.
We present the example of the Thompson group $F=\langle x_0, x_1,
x_2, ... | x_i^{-1} x_n x_i = x_{n+1} $ for all $ i<n$ and $n \in
{\mathbb N}\rangle$, which is of type FP$_\infty$.  It has periodic
cohomology with integral coefficients and vanishing cohomology with
${\mathbb Z}F$-module coefficients, but it does not have jump
cohomology over any ring $R$.

Next we apply Theorem 1.5 to groups that can act freely and properly
on a finite dimension complex homotopy equivalent a closed,
orientable manifold.

\vspace{.2cm}

\begin{cor}\label{LHS} If $G$ is a solvable group, that acts freely and
properly on an $(m+n)$-dimensional complex homotopy equivalent to
closed, orientable manifold of dimension n, then the Hirsch rank
$h(G)\leq m$.\end{cor}

\begin{cor}\label{LHS} If $G$ has a non-countable, torsion-free, solvable
subgroup, then $G$ cannot act freely and properly on a finite
dimensional complex homotopy equivalent to a closed, orientable
manifold.\end{cor}

\vspace{.2cm}

Lastly, in section 4, we find an obstruction for a countable group
acting on certain complexes to be free. We apply the next
proposition to groups that admit free and proper actions on finite
dimension complexes homotopy equivalent to closed, orientable
manifolds.

\vspace{.2cm}

\begin{prop}\label{LHS} Let $G$ be a finitely generated torsion-free group that acts
freely and properly on an $n$-dimensional CW-complex $X$ such that
the $G$-invariant submodule $(H_{n-1}(X,R))^G\ne 0$ . Then $G$ is a
free group if and only if Ext$_{RG}^2(R,H_{n}(X,R))$ is trivial.
\end{prop}

\vspace{0.2cm}

\begin{thm}\label{LHS} If $G$ is a torsion-free countable group that acts freely and
properly on an $n$-dimensional complex homotopy equivalent to a
closed, orientable manifold of dimension $n-1$, then $G$ must be
free. \vspace{0.3cm}
\end{thm}

I would like to thank Alejandro Adem, Donald Passman and Olympia
Talelli for their conversations and advice. Alejandro Adem and
Olympia Talelli read an earlier version of the paper and made useful
suggestions.

\section{Periodicity and Boundedness}

\vspace{0.3cm}

Unless otherwise specified, $R$ will denote any commutative ring
with a unit.

A group $G$ is said to be $R$-{\it torsion-free} if for every finite
subgroup $H$ of $G$, $n\cdot 1_R$ is an invertible element in $R$,
where $n$ is the order of $H$. One should observe that this is
equivalent to the following definition, which will be used in Lemma
3.4.

$G$ is $R$-torsion-free if and only if for each prime $p$, such that
there exists a non-identity element $g\in G$ with $g^p=1$, $p\cdot
1_R$ is an invertible element in $R$.

\vspace{0.2cm}
\begin{defn}\label{LHS} Cohomological dimension of a discrete group $G$ over $R$, denoted
$cd_R(G)$, is defined as

\vspace{0cm}

\begin{center}
inf$\{n:$ Ext$_{RG}^i(R,-)=0$ for $i>n\}$
\end{center}

\noindent and the homological dimension, $hd_R(G)$, of $G$ over $R$
is

\begin{center}
inf$\{n:$ Tor$^{RG}_i(R,-)=0$ for $i>n\}$
\end{center}
\end{defn}

An immediate consequence of $G$ having finite cohomological or
homological dimension is that $G$ must be $R$-torsion-free.

\vspace{0.2cm}

\begin{defn}\label{LHS} A discrete group $G$ has {\it periodic cohomology over} $R$ if
there is a cohomology class $\alpha\in$ Ext$_{RG}^{\ast}(R,R)$ with
$|\alpha|>0$ and an integer $n\geq 0$ such that the cup product map
\begin{center}
 $\alpha \cup -:$ Ext$_{RG}^{i}(R, M)\rightarrow$ Ext$_{RG}^{i+|\alpha|}(R, M)$
\end{center}

\noindent is an isomorphism for every $RG$-module $M$ and every
integer $i\geq n$. {\normalfont The smallest degree of such $\alpha$
will be the ${\it period}$ of the periodic cohomology of $G$.
If $R={\mathbb Z}$, we simply say that $G$ has {\it
periodic cohomology}}.\end{defn}

The above notion of periodicity was used by Adem and Smith in
\cite{AS}. They prove a generalization of Wall's conjecture
\cite{MT} assuming that the periodicity is induced by the cup
product with a cohomology class in $H^*(G,\mathbb Z)$.

\vspace{0.1cm}
\begin{thm}\label{LHS} {\normalfont (Adem-Smith)}. A discrete group $G$ has
periodic cohomology if and only if $G$ acts freely and properly on a
finite dimensional complex homotopy equivalent to a sphere.
\end{thm}

\vspace{0.2cm}

Let us investigate a weaker cohomological condition, which arises
naturally, when a group $G$ has a periodic cohomology over a ring
$R$. It will be our main interest in section 3, where we will
investigate specific types of groups satisfying the condition.
\vspace{0.2cm}

\begin{defn}\label{LHS} A discrete group $G$ has {\it jump cohomology (homology) over} $R$ if there
exists an integer $k\geq 0$, such that if $H$ is any subgroup of $G$
with $cd_R(H)<\infty$ ($hd_R(H)<\infty$), then $cd_R(H)\leq k$
($hd_R(H)\leq k$). {\normalfont The bound $k$ will be called the
{\it jump height over R}. If $R={\mathbb Z}$ then we simply say $G$
has {\it jump cohomology (homology)} and a {\it jump height k}.}
\end{defn}

This next proposition gives a more geometric meaning to the jump height.
\vspace{0.2cm}

\begin{prop}\label{LHS}
 Let $G$ be a discrete group acting freely and properly on an
$e$-dimensional CW-complex $E$. Suppose there is an integer $n\leq e$
such that $H_i(E,{\mathbb Z})=0$ for all $i>n$, and
$H_n(E,{\mathbb Z})\cong {\mathbb Z}$. Set $k=e-n$.
Then $G$ has jump cohomology and homology of height $k$.
\end{prop}

\noindent
\begin{proof} We will show the proof of the proposition concerning cohomology.
The proof of the homological version is similar.

Suppose $H$ is a subgroup of $G$ with finite cohomological dimension
$m$. We can assume that the induced action of $H$ on
$H_n(E,\mathbb{Z})$ is trivial, otherwise we can pass to an index
two subgroup of $H$ and use Serre's Theorem, \cite{ser}.

We get the following fibration
\begin{center}
$E\longrightarrow E\times_{H}EH{\buildrel\pi\over\longrightarrow}
BH$
\end{center}

\noindent where $E\times_{H}EH$ is the Borel Construction and it is
homotopy equivalent to $E/H$. Therefore, by the Leray-Serre spectral
sequence we have
$$E^{p,q}_2=H^p(BH, H^q(E, M))\Longrightarrow H^{p+q}(E/H,M)$$

There exists a ${\mathbb Z}G$-module $F$ such that
$H^{m}(H,F)\ne 0$. By the usual corner argument of the sequence it
follows
$$H^{m+n}(E/H, F)\cong H^m(BH, H^n(E, F))$$
By the Universal Coefficient Theorem we have the following
isomorphism of ${\mathbb Z}H$-modules
\begin{center}
$H^n(E,F)\cong$ Hom$(H_n(E, {\mathbb Z}),F)\oplus N \cong F \oplus N$
\end{center}
\noindent
where $N=$ Ext$(H_{n-1}(E, {\mathbb Z}), F)$. This shows that
$$H^{m+n}(E/H, F)\cong H^m(BH, F)\oplus H^m(BH, N)$$
Therefore, $H^{m+n}(E/H, F)\ne 0$ and $m+n\leq e$.

\end{proof}

\vspace{0.1cm}
\begin{remk}\label{LHS}{\normalfont If there exist a finite dimensional
free-$G$-CW-complex $E$ homotopy equivalent to a sphere $S^n$, then
$G$ has a periodic cohomology over any ring $R$. For we can always
construct the join $E\ast E$ with the diagonal $G$-action. This
complex is homotopy equivalent to $S^{2n+1}$ and $H_*(E\ast
E,\mathbb Z)$ has the trivial $G$-action. Therefore the induced
$G$-action on $H^*(E\ast E,R)$ is also trivial. Then the Gysin exact
sequence of cohomology
over $R$ shows that the cup product with the Euler class of the
spherical fibration induces an isomorphism in cohomology of $G$
after dimension $2k$, where $2k$ is the jump height over $R$. This,
together with Theorem 2.3, yields that, if a group $G$ has periodic
cohomology, then it has periodic cohomology over any ring $R$.}
\end{remk}

Now, suppose $G$ has periodic cohomology over $R$ starting in
dimension $k+1$. If $H$ is any subgroup of $G$ with finite
cohomological dimension $h$ over $R$, using Shapiro's Lemma, we have

\begin{center}
Ext$_{RH}^{k+1}(R,M)\cong$ Ext$_{RG}^{k+1}(R, $ Coind$_H^G M)
\displaystyle{{\buildrel{\alpha^h \cup - } \over \longrightarrow}}$

Ext$_{RG}^{k+h|\alpha|+1}(R, $ Coind$_H^G M)\cong$
Ext$_{RH}^{k+h|\alpha|+1}(R,M)=0$
\end{center}

\noindent for any $RH$-module $M$. This establishes the following
lemma. \vspace{0.1cm}
\begin{lem}\label{LHS} If a group $G$ has periodic cohomology over $R$ starting in
dimension $k+1$, then $G$ has a jump cohomology of height $k$ over
$R$.\end{lem}
\qed \vspace{0.1cm}

Therefore, among the three cohomological conditions stated in the
conjecture, jump cohomology is the weakest.

\section{Some classes of groups}

\vspace{0.3cm}

\noindent The following class of groups was introduced by P.
Kropholler in \cite{krop}. Our next theorem shows that they satisfy
the conjecture stated in section 1.

\vspace{.3cm}

\begin{defn}\label{LHS} Let $\mathcal{X}$ denote a class of groups. Define ${\sl
H}\mathcal{X}$ to be the smallest class of groups containing
$\mathcal{X}$ with the property: if a group $G$ acts cellularly on a
finite dimensional contractible CW-complex with all isotropy
subgroups in ${\sl H}\mathcal{X}$, then $G$ is in ${\sl
H}\mathcal{X}$.\end{defn}

In \cite{krop} many properties of these classes of groups, such as
subgroup and extension closure, closure under countable direct
unions and free product, were shown. The main interest to us is the
hierarchical description of ${\sl H}\mathcal{X}$-groups defined by
operations ${\sl H}_{\alpha}$ for each ordinal $\alpha$ inductively:

\noindent -${\sl H}_0 {\mathcal X}= {\mathcal X}$

\noindent For each $\beta > 0$,

\noindent -${\sl H}_{\beta} \mathcal{X}$ is the class of groups $G$
which act cellularly  on a finite dimensional contractible
CW-complex $X$ such that for each cell $\sigma$ of $X$ the isotropy
group $G_{\sigma}$ is in ${\sl H}_{\alpha} \mathcal{X}$ for some
$\alpha < \beta$.

\vspace{.2cm}

It is immediate that a group $G$ is a ${\sl H}\mathcal{X}$-group if
and only if there is an $\alpha$ such that $G$ is in ${\sl
H}_{\alpha} \mathcal{X}$. Let $\mathcal{F}$ denote the class of all
finite groups. We obtain the following result pertaining to ${\sl
H}\mathcal{F}$-groups.

\vspace{0.3cm}

\begin{thm}\label{LHS} Let $G$ be a discrete group without $R$-torsion with jump
cohomology of height $k$ over $R$. If $G$ is an ${\sl
H}\mathcal{F}$-group, then $cd_R(G)\leq k$. In particular, any ${\sl
H}\mathcal{F}$-group $G$ has jump cohomology of height $k$ over
$\mathbb Q$ if and only if $cd_{\mathbb Q}(G)\leq k$. \end{thm}

First we need a lemma. \vspace{0.1cm}
\begin{lem}\label{LHS} If a discrete group $G$ acts cellularly on an
$n$-dimensional contractible complex $X$, such that for any cell
stabilizers of $X$, $cd_R(G_{\sigma})\leq k$, then $cd_R(G)\leq
k+n$.\end{lem}

\noindent
\begin{proof} By considering the double complex $Hom_{RG}(P_*,C^*_R(X,M))$,
where $P_*$ is a projective resolution of $R$ over $RG$ and
$C^*_R(X,M)$ is the cellular cochain complex of $X$ with $RG$-module
coefficients $M$ (see for example \cite{brown}) we can derive the
first-quadrant spectral sequence

\begin{center}
$E^{p,q}_1= $ Ext$^q_{RG}(R, C^p_R(X,M))\Longrightarrow $
Ext$^{p+q}_{RG}(R,C^*_R(X,M))$
\end{center}

\noindent where Ext$^*_{RG}(R,C^*_R(X,M))$ is the cohomology of the
total complex associated to the double complex.

For each $p$-cell $\sigma$ of $X$ there exists a $G_{\sigma}$-module
$R_{\sigma}$. This module is isomorphic to $R$ additively, and
$G_{\sigma}$ acts on it through the orientation character. Let
$M^{\sigma}=Hom_R(R_{\sigma}, M)$. Let $X_p$ denote the collection
of all the $p$-cells and let $\Sigma_p$ be a set of representatives
of all the $G$-orbits in $X_p$. We have the following decomposition

\begin{center}
$C^p_R(X,M)=Hom_R(C_p^R(X),M)=\displaystyle{\bigoplus_{\sigma \in
X_p}}Hom_R(R_{\sigma},M)=$

$\displaystyle{\bigoplus_{\sigma \in
{\Sigma}_p}}Hom_R($Ind$^G_{G_{\sigma}}R_{\sigma},M)=\displaystyle{\bigoplus_{\sigma
\in {\Sigma}_p}}$Coind$^G_{G_{\sigma}} M^{\sigma}$
\end{center}

Now, by Shapiro's lemma,

\begin{center}
Ext$^q_{RG}(R, C^p_R(X,M))\cong \displaystyle{\bigoplus_{\sigma \in
\Sigma_p}} $Ext$^q_{R[G_{\sigma}]}(R, M^{\sigma})$.
\end{center}

Since $X$ is contractible  the two cochain complexes $C^*_R(X,M)$
and $C^*_R(pt.,M)$ are homotopy equivalent. Hence, \vspace{-.2cm}

\begin{center}
Ext$^*_{RG}(R,C^*_R(X,M))\cong$ Ext$^*_{RG}(R,C^*_R(pt.,M))= $
Ext$^*_{RG}(R,M)$.
\end{center}

The spectral sequence then becomes
\begin{center}
$E^{p,q}_1=\displaystyle{\bigoplus_{\sigma \in \Sigma_p}}
$Ext$^q_{R[G_{\sigma}]}(R, M^{\sigma})\Longrightarrow $
Ext$^{p+q}_{RG}(R,M)$
\end{center}

It follows that $E^{p,q}_1 = 0$ if $p>n$ or $q>k$. We can infer
inductively for $i=0,1,2,...$ that the differentials $d^{p,q}_i = 0$
and the terms of the spectral sequence $E^{p,q}_i = 0$ if $p>n$ or
$q>k$. Therefore, $E^{p,q}_{\infty} \cong E^{p,q}_1 = 0$ if $p>n$ or
$q>k$ and Ext$^{p+q}_{RG}(R,M)=0$ if $p>n$ and $q>k$. Thus
$cd_R(G)\leq n+k$.

\end{proof}

\vspace{.2cm}

\noindent
\begin{proof}[Proof of Theorem] We will show that $cd_R(G)\leq k$
by transfinite induction, using the hierarchical description of
${\sl H}\mathcal{F}$-groups.

If $G$ is in ${\sl H}_0 \mathcal{F} = \mathcal{F}$, then $G$ is a
finite group. So, $cd_R(G)=0$. Assume now that for a fixed $\beta$,
and all $\alpha < \beta$, $cd_R(H) \leq k$ for any subgroup $H$ of
$G$ in ${\sl H}_{\alpha}\mathcal{F}$. Suppose $G$ is an ${\sl
H}_{\beta}\mathcal{F}$-group. By assumption there is a finite
dimensional contractible $G$-CW-complex $X$ such that for any cell
$\sigma$ of $X$ the isotropy group $G_{\sigma}$ is an ${\sl
H}_{\alpha}\mathcal{F}$-group for some $\alpha < \beta$. So, by
induction, $cd(G_{\sigma})\leq k$. Let $n$ be the dimension of $X$.
Now, by the lemma, $cd_R(G)\leq k+n$. Hence, by hypothesis,
$cd_R(G)\leq k$.

\end{proof}

\vspace{0.2cm}

Our next goal is to study solvable groups with jump cohomology. For
this we need to assume $R$ is an integral domain. Note that by the
previous theorem the conjecture stated in section 1 is true for
countable $R$-torsion-free solvable groups, for they are ${\sl
H}\mathcal{F}$-groups, for an arbitrary commutative ring $R$. What
proceeds extends this to uncountable solvable groups when $R$ is a
domain. In fact, our arguments show that if a torsion-free solvable
group has jump cohomology, then it is countable, since it must have
finite Hirsch rank. We denote by ${\mathbb F}$ the fraction field of
$R$, so ${\mathbb F}$ is a field of an arbitrary characteristic.
\vspace{0.2cm}
\begin{lem}\label{LHS} Let $1=G_0 \subset G_1 \subset $ ... $\subset G_n=G$
be an upper central series of $G$. If $G$ is an $R$-torsion-free
group, then so are the factor groups $G_i/G_{i-1}$,
$i=1,..,n$.\end{lem}

\vspace{0.2cm}

\noindent
\begin{proof}Recall that an upper central series for $G$ can be defined
inductively by letting $G_i$ to be the normal subgroup of $G$ such
that $G_i/G_{i-1} = Z(G/G_{i-1})$, where $Z(G/G_{i-1})$ is the
center of
 the quotient group $G/G_{i-1}$.

Let us proceed by induction. Suppose, for all $j\leq i$,
$G_j/G_{j-1}$ are $R$-torsion-free. We need to show that the same is
true for the quotient group $G_{i+1}/G_i$. Suppose $\bar x \in
G_{i+1}/G_i$ such that $\bar x ^p=1$ for a prime number $p$, and
$p\cdot 1_R$ is not invertible in the ring $R$. This shows that
$x^p\in G_i$. Now for any element $g$ of $G$, $[g,x]^p\equiv
[g,x^p]\equiv 1$ mod $G_{i-1}$. Since $[g,x] \in G_i$, the induction
shows $[g,x] \in G_{i-1}$. So, in the quotient group $G/G_{i-1}$,
$[\bar g, \bar x]=1$ for all $\bar g \in G/G_{i-1}$. Therefore, by
definition, $\bar x \in G_i/G_{i-1}$ and $x \in G_i$. It follows
that $\bar x$ is the identity element in $G_{i+1}/G_i$.

\end{proof}

\vspace{0.2cm}
\begin{prop}\label{LHS} Let $G$ be a nilpotent group without $R$-torsion
with jump homology of height $k$ over $R$. Then $hd_R(G)=h(G)\leq k$
and Tor$^{RG}_{h(G)}(R,{\mathbb F})\cong {\mathbb F}$, where $h(G)$
is the Hirsch rank of $G$.\end{prop}

\begin{proof}  Let $1=G_0 \subset G_1 \subset $ ... $\subset G_n=G$
be a central series of $G$ with $R$-torsion-free factor groups. We use induction on $n$.

Let us assume $h(G_i)= hd_R(G_i)=h_i$,  and
Tor$^{RG_i}_{h_i}(R,{\mathbb F})\cong{\mathbb F}^{\omega_i}$, where
${\mathbb F}^{\omega_i}$ is additively isomorphic to ${\mathbb F}$
with a trivial $G/G_i$-action. Let $A\subset G_{i+1}/G_i$ and
$A\cong \mathbb{Z}^m$. Let ${{G_{i+1}^\prime}}=\varphi^{-1}(A)$,
where $\varphi:G_{i+1}\rightarrow G_{i+1}/G_i$ is the canonical
quotient map. So, there exists an extension, $1\rightarrow G_i
\rightarrow {{G_{i+1}^\prime}}\rightarrow A\rightarrow 1$.

By the corner argument of the Lyndon-Hochschild-Serre spectral
sequence of the extension above it follows,
\begin{center}
Tor$^{RG_{i+1}^\prime}_{m+h_i}(R,{\mathbb F})\cong$
Tor$^{RA}_{m}(R,$ Tor$^{RG_i}_{h_i}(R,{\mathbb F}))\cong$
Tor$^{RA}_m(R, {\mathbb F}^{\omega_i})$
\end{center}

Since $A\cong {\mathbb Z}^m$,
\begin{center}
Tor$^{RA}_m(R,{\mathbb F}^{\omega_i})\cong
\bigwedge^m (A\otimes {\mathbb F}^{\omega_i})\cong \bigwedge^m
(\mathbb{Z}^m \otimes {\mathbb F}^{\omega_i})\cong {\mathbb F}$
\end{center}
\noindent where $\bigwedge^\ast (A\otimes {\mathbb F}^{\omega_i})$
is the exterior algebra of the ${\mathbb F}$-vector space $A\otimes
{\mathbb F}^{\omega_i}$. This shows
Tor$^{RG_{i+1}^\prime}_{m+h_i}(R,{\mathbb F})\ne 0$. On the other
hand, since $hd_R(G_i)=h_i$ and $hd_R(A)=m$, the corner argument of
the spectral sequence also shows, $hd_R({G_{i+1}^\prime})\leq m +
h_i$. Thus $hd_R({G_{i+1}^\prime})= m+h_i$. So, by the hypothesis
$m+h_i\leq k$, showing $m\leq k$. This proves that $G_{i+1}/G_i$ is
a finite-rank abelian group. Thus, $G_{i+1}$ has finite Hirsch rank,
and since $G_{i+1}/G_i$ is also $R$-torsion-free, $h(G_{i+1}/G_i)=
hd_R(G_{i+1}/G_i)$. So, again by the corner argument of the spectral
sequence associated to the extension, $1\rightarrow G_i \rightarrow
{G_{i+1}}\rightarrow G_{i+1}/G_i\rightarrow 1$, we have
$hd_R(G_{i+1})\leq h(G_{i+1})$.

Let $m=hd_R(G_{i+1}/G_i)$. The spectral sequence gives

\begin{center}
Tor$^{RG_{i+1}}_{m+h_i}(R,{\mathbb F})\cong$
Tor$^{R[G_{i+1}/G_i]}_{m}(R, $ Tor$^{RG_i}_{h_i}(R, {\mathbb
F}))\cong$ Tor$^{R[G_{i+1}/G_i]}_m(R,{\mathbb F}^{\omega_i})\cong
\bigwedge^m (G_{i+1}/G_i\otimes {\mathbb F}^{\omega_i})$
\end{center}

The last isomorphism follows from the fact that $G_{i+1}/G_i$ is an
$R$-torsion-free abelian group and $\mathbb F^{\omega_i}$ has
trivial $G_{i+1}/G_i$-action. Also, since $G_{i+1}/G_i$ is central
in $G/G_i$, the induced action of $G/G_{i+1}$ on $\bigwedge^m
(G_{i+1}/G_i\otimes {\mathbb F}^{\omega_i})$ is trivial. So the
module $\bigwedge^m (G_{i+1}/G_i\otimes {\mathbb F}^{\omega_i})$ is
${\mathbb F}^{\omega_{i+1}}$, and $hd_R(G_{i+1})=h(G_{i+1})$. This
completes the induction.

\end{proof}

Our next result requires a version of a theorem of Stammbach
\cite{stam}. For the reader's convenience we outline its proof which
can be also found in \cite{bieri}.

\begin{prop}\label{LHS}{\normalfont (Stammbach)} Let $G$ be an
$R$-torsion-free solvable group with the Hirsch rank $h$, where $R$
is an integral domain of characteristic zero. Then $hd_R(G)=h$ and
Tor$^{\mathbb FG}_h(\mathbb F, A)\cong {\mathbb F}$, for some
$\mathbb FG$-module $A$ isomorphic to ${\mathbb F}$ as an $\mathbb
F$-module.\end{prop}

\noindent
\begin{proof}Let $1=G_0 \subset G_1 \subset $ ... $\subset G_n=G$
be a derived series of $G$. Let $S_i = G_i/G_{i-1}$ and
$rk(S_i)=h_i$. Set $F_i = $ Tor$^{\mathbb FS_i}_{h_i}(\mathbb F,
\mathbb F)$. From Proposition 3.5 it follows that $F_i$ is
additively isomorphic to ${\mathbb F}$.

Denote by $F_i^{op}$ the additive group of $F_i$ with the inverse
$G$-action. Let $A= \otimes_{i=1}^{i=n}F_i^{op}$. We use induction
on $n$ to prove Tor$^{\mathbb FG}_h(\mathbb F, A)\cong {\mathbb F}$.
For $n=1$ the result is clear. Suppose $n\geq 2$, then $G/G_1$ has
derived length $n-1$. Thus, by the Lyndon-Hochschild-Serre spectral
sequence of the extension $1\to G_{1}\to G\to G/G_{1}\to 1$, we have

\begin{center}
Tor$_h^{\mathbb FG}(\mathbb F, A)\cong$ Tor$_{h-h_1}^{\mathbb
F(G/G_1)}(\mathbb F,$ Tor$_{h_1}^{\mathbb FG_1}(\mathbb F, A))\cong$

Tor$_{h-h_1}^{\mathbb F(G/G_1)}(\mathbb F,$ Tor$_{h_1}^{\mathbb
FG_1}(\mathbb F, {\mathbb F})\otimes A) \cong $
Tor$_{h-h_1}^{\mathbb F(G/G_1)}(\mathbb F,
\otimes_{i=2}^{i=n}F_i^{op})\cong {\mathbb F}$. \end{center}

Now the assertion that $hd_R(G)=h$ follows from the fact that $h\leq
hd_{\mathbb F}(G)\leq hd_R(G)\leq h$, when $G$ is $R$-torsion-free.
\end{proof}
\begin{thm}\label{LHS} Let $G$ be a solvable group which acts freely
and properly
on an $e$-dimensional CW-complex $E$. Suppose there is an integer $n\leq e$
such that $H_i(E,{\mathbb Z})=0$ for all $i>n$, and $H_n(E,{\mathbb Z})\cong
{\mathbb Z}$, then $h(G)\leq e - n$. If $G$ is also without
$R$-torsion, $hd_R(G)\leq e - n$.\end{thm}
\begin{proof} Suppose $h(G)> e - n$. Then, there exists a subgroup $H$ of
$G$ which has finite Hirsch rank larger than $e-n$. Now by 3.6 it
follows that $hd_{\mathbb Q}(H)= h(H)$. Therefore the homological
dimension of $H$ over $\mathbb Q$ is finite and larger than $e-n$,
which contradicts Proposition 2.5.
\end{proof}

\begin{examp}\label{LHS} {\normalfont Let $F$ denote the Thompson group
 defined by the presentation

\begin{center}
$\langle x_0, x_1, x_2, ... | x_i^{-1} x_n x_i = x_{n+1} $ for all $
i<n$ and $n\in {\mathbb N}\rangle$
\end{center}

This group is finitely presented for it also has the presentation

\begin{center}
$\langle x_0, x_1 | x_2^{x_0}=x_3, x_3^{x_1}=x_4\rangle$
\end{center}

\noindent Brown and Goeghegan \cite{BG} showed that it is also of
type FP$_\infty$, i.e. it has a projective ${\mathbb Z}G$-resolution
of ${\mathbb Z}$, where each module is finitely generated. They also
proved $H_n(F,{\mathbb Z})={\mathbb Z}\oplus {\mathbb Z}$ for $n\geq
1$ and $H^n(F,{\mathbb Z}F)=0$ for all $n$. This may suggest that
this group has periodic cohomology. But in fact, $F$ cannot have
periodic or jump cohomology over any ring, since it has an infinite
rank abelian subgroup ${\mathbb Z}^\infty$. Namely, $\langle x_0
x_1^{-1}, x_2 x_3^{-1},  x_4 x_5^{-1}, ... \rangle \subseteq
F$.}\end{examp}

This reinstates our belief that periodicity and jumps are rather
special conditions on the cohomology of a discrete group, since
countable groups containing the Thompson group are the only class of
groups, known to us, not to be ${\sl H}\mathcal{F}$.

\vspace{0.2cm}
\begin{cor}\label{LHS} If $G$ is a solvable group that acts freely and properly
on an $(m+n)$-dimensional complex homotopy equivalent to a closed,
orientable manifold of dimension n, then the Hirsch rank $h(G)\leq
m$.\end{cor} \qed
\begin{cor}\label{LHS} If $G$ has a non-countable, torsion-free, solvable
subgroup, then $G$ cannot act freely and properly on a finite
dimensional complex homotopy equivalent to a closed, orientable
manifold.\end{cor}

\begin{proof} Let $H$ be a non-countable, torsion-free, solvable subgroup of
$G$. By \cite{bieri}, Lemma 7.9, p. 100, $H$ has infinite Hirsch
rank. Therefore it cannot act freely and properly on a finite
dimensional complex homotopy equivalent to a closed, orientable
manifold.
\end{proof}

\section{Actions of Height One}

\vspace{0.3cm}

These next results apply only to countable torsion-free groups. We
find an obstruction for a group to be free in a specific setting. It
would be interesting to know whether the following proposition holds
for not necessarily countable groups and what implications there are
when the obstruction does not vanish.

\vspace{0.3cm}
\begin{prop}\label{LHS}  Let $G$ be a finitely generated torsion-free group that acts
freely and properly on an $n$-dimensional CW-complex $X$ such that
the $G$-invariant submodule $(H_{n-1}(X,R))^G\ne 0$. Then $G$ is a
free group if and only if Ext$_{RG}^2(R,H_{n}(X,R))=0$.\end{prop}

\begin{proof} If $G$ is a free group then $cd_R(G)\leq 1$, so  Ext$_{RG}^2(R,H_{n}(X,R))$
must be trivial.

Conversely, suppose Ext$_{RG}^2(R,H_{n}(X,R))=0$. Let $(C_*(X),d_*)$
be the chain complex of $X$, where $C_i(X)$ is the free $R$-module
generated by the $i$-cells of $X$.  We proceed by induction on the
number of generators of $G$ with the case of $G$ having exactly one
generator understood.

There exists an exact sequence of ${\mathbb Z}G$-modules

\begin{center}
$0\rightarrow d_n(C_n(X))\rightarrow $ker$(d_{n-1})\rightarrow
H_{n-1}(X,R)\rightarrow 0$
\end{center}

\noindent and an associated long exact sequence

\begin{center}
Ext$_{RG}^0(R,$ ker$(d_{n-1}))\rightarrow$ Ext$_{RG}^0(R,
H_{n-1}(X,R))\rightarrow$ Ext$_{RG}^1(R, d_n(C_n(X)))\rightarrow$
Ext$_{RG}^1(R,$ ker$(d_{n-1}))$
\end{center}

The first term of the long exact sequence is the $G$-invariant
submodule (ker$(d_{n-1}))^G$. This module is trivial, for it is a
submodule of the $G$-invariant submodule $(C_{n-1}(X))^G$. But $G$
is infinite and $C_{n-1}(X)$ is a free $RG$-module, hence
$(C_{n-1}(X))^G=0$. Thus the second term, the $G$-invariant
submodule $(H_{n-1}(X,R))^G$, of the long exact sequence injects
into Ext$_{RG}^1(R,d_n(C_n(X)))$. This shows that
Ext$_{RG}^1(R,d_n(C_n(X)))\ne 0$.

On the other hand there exists an exact sequence of $RG$-modules

\begin{center}
$0\rightarrow H_n(X,R)\rightarrow {C_n(X)}\rightarrow
d_n(C_n(X))\rightarrow 0$
\end{center}

\noindent The long exact sequence associated to this gives:

\begin{center}
Ext$_{RG}^1(R, H_n(X,R))\rightarrow$
Ext$_{RG}^1(R,{C_n(X)})\rightarrow$
Ext$_{RG}^1(R,d_n(C_n(X)))\rightarrow$ Ext$_{RG}^2(R, H_n(X,R))$
\end{center}

\noindent According to our assumption Ext$_{RG}^2(R, H_n(X,R))=0$,
so the term Ext$_{RG}^1(R,{C_n(X)})$ surjects onto the module
Ext$_{RG}^1(R,d_n(C_n(X)))$, which is nontrivial. Therefore
Ext$_{RG}^1(R,{C_n(X)})\ne 0$. Now since $G$ is finitely generated
and $C_n(X)$ is a free $RG$-module, we must have
Ext$_{RG}^1(R,RG)\ne 0$. So by \cite{cohen} $G$ must split as a
nontrivial free product of groups.

Suppose then $G=H\ast K$ where $H,K\ne G$. For any $RG$-module $M$,

\begin{center}
Ext$_{RG}^2(R,M)=$ Ext$_{RH}^2(R,M)\oplus$ Ext$_{RK}^2(R,M)$
\end{center}

\noindent It then follows that,

\begin{center}
Ext$_{RH}^2(R, H_n(X,R))=$ Ext$_{RK}^2(R, H_n(X,R))=0$.
\end{center}

Also by Grushko's Theorem \cite{cohen} there exist a fewer number of
generators for  each of the groups $H$ and $K$ than the number of
generators of $G$. So we can apply induction on these subgroups to
show that they are free.

\end{proof}

\begin{thm}\label{LHS} If $G$ is a torsion-free countable group that acts
freely and properly on an $n$-dimensional complex homotopy
equivalent to a closed, orientable manifold of dimension $n-1$, then
$G$ must be free.\end{thm}

\begin{proof} First, let us assume $G$ is finitely generated.
If the induced action of $G$ on $H_{n-1}(X)= \mathbb Z$ is
nontrivial then $G$ must contain an index two subgroup $H$ which
acts trivially on $H_{n-1}(X)$. Thus, there always exists such
$H\subseteq G$ acting trivially, with $\lbrack G:H \rbrack \leq 2$.
Also, since $H_n(X)=0$, $H^2(H,H_n (X))=0$. By Proposition 4.1 it
follows that $H$ is free.
 Since $H$ is a finite index subgroup, $G$ must be free.

In general, since $G$ is a countable group, it is a countable union
of an ascending chain of finitely generated subgroups. By the
Theorem of Berstein, \cite{grun}, we have $cd(G)\leq 2$. Using
Proposition 2.5, it follows that $cd(G)\leq 1$. Therefore $G$ is
free.

\end{proof}

\end{document}